\documentclass[12pt]{article}
\usepackage[english]{babel}
\usepackage{amsmath,amsthm}
\usepackage{amsfonts}

\theoremstyle{definition}

\theoremstyle{remark}

\numberwithin{equation}{section}

\begin{document}

\begin{center}
\bfseries  Riemann hypothesis is not correct   \\
\end{center}

\begin{center}

\vspace{5mm}
JINHUA  FEI\\
\vspace{5mm}

ChangLing  University   \, Baoji \,  Shannxi  \,  P.R.China \\
\vspace{5mm}

E-mail: feijinhuayoujian@msn.com \\

\end{center}

\vspace{3mm}

{\bfseries Abstract.}\,  This paper use   Nevanlinna's Second Main Theorem  of  the  value distribution theory,   we got an important conclusion  by  Riemann hypothesis. this conclusion contradicts the Theorem 8.12 in Titchmarsh's book "Theory of the Riemann Zeta-functions", therefore we prove that Riemann hypothesis is incorrect.
\\


{\bfseries MR(2000)  Subject  Classification \quad 30D35,  11M06 } \\

\vspace{8mm}

    First, we give some notations, definitions and theorems in the theory of value distribution, its contents see the references [1] and [2]. \\

 We write

\[  \log^+ x  =   \left\{
     \begin{array}{ccc}
     \log x  \qquad  \quad  1 \leq  x \\
     \quad 0    \qquad  \quad    0 \leq x < 1
     \end{array}
     \right.
\]

It is easy to see that  $  \,\,  \log x  \, \leq  \, \log^+ x   $.\\

  Let $ f(z)  $ is a non-constant  meromorphic function  in the  circle $ \,\, |z|  \, <  \, R \, ,  \,\,  0 \,  < R  \,  \leq  \, \infty  \, $.
  $ n(r,f) $  represents  the number of  poles of  $ f(z) $ on  the  circle  $ |z| \, \leq  \, r \,( \, 0 <  r < R \, ) \,  $,  the multiplicity of poles is included. $ n(0,f) $  represents the order of pole of $ f(z)  $  in the  origin.  For arbitrary complex number  $ a  \neq  \infty  \, , \,\,\,\,  n( r , \frac{1}{ f-a  }) \,\, $ represents  the number of  zeros of  $ f(z) - a $ in the circle $ \, |z| \leq  r \,\, (  \, 0 <   r  <  R \, ) \, $, the multiplicity of zeros  is included. $ \,\,  n( 0 , \frac{1}{f-a} ) \,\, $   represents the order of zero of $ \,\, f(z) - a \,\, $ in the  origin.\\

We write
$$  m( r, f)  =  \frac{1}{ 2 \pi}  \, \int_0^{2 \pi} \, \log^+  \left| f( r e^{i \varphi} )  \, \right|  d \varphi  $$

$$ N ( r , f )\, = \, \int_0^r \, \frac{ n(t,f) - n(0,f)}{ t } \, d t   \,\, + \,\, n(0,f) \, \log \, r $$ \\

and$ \,\, T( r, f  ) \, = \, m( r , f) \, + \, N ( r , f ) $ .

   $\,\,  T ( r, f ) $  is called the characteristic function  of  $ f ( z )  $. \\

{\bfseries\footnotesize LEMMA 1.}  If  $ f(z)  $  is  a  analytical function in the circle  $ \, |z| < R \,\, (\, 0 < R \leq  \infty \, ) $, we have

$$ T ( r , f ) \, \leq \,  \log^+ \, M ( r , f ) \, \leq \,  \frac{ \rho + r }{ \rho - r } \, \, T ( \rho \, , f)    (\, 0 < r < \rho < R  )   $$

where  $ \,\,  M( r, f ) \, = \, \max_{ |z| = r   }  \, | f(z)  | $

The lemma 1 follows from the  References [1], page 57. \\

{\bfseries\footnotesize LEMMA 2.}   Let $ f(z) $  is a non-constant  meromorphic function  in the circle $ |z| < R \,\, ( \, 0 < R \leq \infty  )\, $.
  $ \,\, a_\lambda \,\, ( \, \lambda = 1,2, ..., h \, )   $ and  $\,\, b_\mu \,\,  (\, \mu = 1,2,..., k\, ) \,  $   are the zeros  and  poles  of    $ f(z)  $  in
   the  circle  $ |z| < \rho  \,\, ( \, 0 < \rho < R \,) $  \,\,respectively,  each zero or pole repeated according to their multiplicity, and $ z = 0  $ is neither zero nor pole  of the function   $ f(z)  $,  then, in  the   circle  $ |z| < \rho  $, we have the following formula

$$  \log \, \left| \,  f(0) \,  \right| \, = \,  \frac{1}{ 2 \pi}  \,\int_0^{ 2 \pi } \,\, \log\, \left| f( \rho e^{ i \varphi} )
   \right| d \varphi  \,\, - \,\, \sum_{ \lambda =1}^h \, \log\, \frac{\rho}{|a_\lambda|} \,\, + \,\,
    \sum_{ \mu =1}^k \, \log\, \frac{\rho}{|b_\mu|} $$
this formula is called Jensen formula.

 The lemma 2 follows from the  References [1], page 48.   \\

{\bfseries\footnotesize LEMMA 3.} Let $ f(z) $ is the meromorphic function  in the circle $ \, |z| \leq  R \, $,  and

 $$ f(0) \,  \neq  \,\,  0 ,\, \, \infty ,\, \, 1, \,\,\,\, \,\, f'(0) \, \neq \, 0 $$

when $ \, 0 < r < R \, $, we have

$$  T ( r , f ) \, <  \, 2 \left\{  N ( R, \frac{1}{f}) \, + \, N ( R , f ) \, + \, N ( R , \frac{1}{ f - 1} ) \right\} $$

$$  + \,\, 4\, \, \log^+ |f(0)| \,\, + \,\, 2 \, \, \log^+ \, \frac{1}{ R |f'(0)| } \,\, + \,\, 24\,\, \log \, \frac{R}{R-r}
\,\, + \,\, 2328 $$

This is  a form  of  Nevanlinna's Second Main Theorem.

 The lemma 3 follows from the  References [1], the theorem 3.1 of the page 75.       \\

  Now, we make some preparations.\\

{\bfseries\footnotesize LEMMA 4.}   if $  f(x)$  is  a function of the nonnegative  degressive,  we have

  $$  \lim_{ N \, \rightarrow  \, \infty} \,\,  \left( \, \, \sum_{n=a}^N \, f(n) \,\, - \,\, \int_a^N  f(x) \, d x  \right) \,\,
  = \,\, \alpha $$

where $  \, 0 \leq  \alpha  \leq  f(a)  \, $. in addition, if  $ \, x \rightarrow  \infty \,  $ ,  $ \, f(x) \rightarrow  0  \, $, we have

$$  \left| \,\,  \sum_{a \leq n \leq \xi }  \, f(n) \,\, - \,\,  \int_a^\xi  f(\nu) \, d \nu   \,\, - \,\, \alpha
 \,\,\right| \,\, \leq \, f( \xi - 1 ) \, , \qquad  (  \, \xi \geq a+1  \,  )  $$

 The lemma 4 follows from the  References [3], the theorem 2 of the page 91.         \\

    Let $ \, s =  \sigma  +  it \,  $ is the  complex number,  when $ \, \sigma  > 1  \,  $,  Riemann  Zeta  function  is

$$  \zeta (s) \, = \, \sum_{n=1}^\infty  \,\, \frac{1}{n^s}   $$

When $ \,  \sigma  >  1 \,  $,  we have

$$  \log \zeta (s) \, =  \, \sum_{n=2}^\infty \, \frac{\Lambda(n)}{n^s  \log n}   $$

where $ \, \Lambda (n) \,  $  is  Mangoldt  function. \\

{\bfseries\footnotesize LEMMA 5.}  If $ t $ is any real number, we have \\

(1) $$ 0.0426  \, \leq \, \left| \, \, \log \zeta ( 4+ it ) \, \, \right|  \, \leq \, 0.0824   $$

(2)$$  \left| \,\,  \zeta ( 4 + it  )   \, - \, 1  \,\, \right|  \, \geq \, 0.0426   $$

(3)$$  0.917 \, \leq \, \left| \,\,  \zeta  ( 4 + it ) \,\, \right|  \, \leq  \,  1.0824 $$

(4)$$  \left|  \,\,  \zeta ' ( 4 + it ) \,\,   \right|   \, \geq \, 0.012  $$ \\

{\bfseries\footnotesize  PROOF.}

(1) $$    \left| \, \, \log \zeta ( 4+ it ) \,\, \right|  \, \leq \, \sum_{n=2}^\infty \,  \frac{\Lambda(n)}{n^4  \log n} \, \leq  \,
 \sum_{n=2}^\infty \, \frac{1}{n^4} \, = \, \frac{\pi^4}{90} \, - \, 1\, \leq \, 0.0824  $$

$$   \left| \, \, \log \zeta ( 4+ it ) \,\, \right| \, \geq \, \frac{1}{2^4} \, - \, \sum_{n=3}^\infty \, \frac{1}{n^4} \, =
 \, 1 \,+\, \frac{2}{2^4} \, - \,  \sum_{n=1}^\infty \, \frac{1}{n^4} \, = \, \frac{9}{8} \,- \, \frac{\pi^4}{90} \, \geq \, 0.0426  $$

(2)

$$    \left| \, \,  \zeta ( 4+ it )  - 1 \,\, \right| \, = \, \left| \,\, \sum_{n=2}^\infty \, \frac{1}{n^{4 + it}} \,\, \right|
 \, \geq \, \frac{1}{2^4} \, - \, \sum_{n=3}^\infty  \, \frac{1}{n^4}  \,$$

$$= \, 1+ \frac{2}{2^4} \, - \,  \sum_{n=1}^\infty  \, \frac{1}{n^4} \, = \, \frac{9}{8} \, - \, \frac{\pi^4}{90} \, \geq \, 0.0426   $$

(3)

$$   \left| \, \, \zeta ( 4+ it ) \,\, \right| \, = \, \left| \,\sum_{n=1}^\infty  \, \frac{1}{n^{4+it}} \, \right| \, \leq \,
 \sum_{n=1}^\infty  \, \frac{1}{n^4} \, = \, \frac{\pi^4}{90} \, \leq \, 1.0824  $$

$$  \left| \, \, \zeta ( 4+ it ) \,\, \right| \, = \, \left| \,\sum_{n=1}^\infty  \, \frac{1}{n^{4+it}} \, \right| \, \geq  \, 1 \,
 - \,  \sum_{n=2}^\infty  \, \frac{1}{n^4} \, = \, 2 \, - \, \sum_{n=1}^\infty  \, \frac{1}{n^4} \, = \, 2\, - \, \frac{\pi^4}{90} \,
  \geq \, 0.917  $$

(4)

$$    \left| \, \, \zeta '( 4+ it ) \,\, \right| \, =  \, \left| \,\, \sum_{n=2}^\infty \, \frac{\log n }{n^{4+it}}
\,\,\right| \, \geq \, \frac{\log2}{2^4} \, - \,  \sum_{n=3}^\infty \, \frac{\log n }{n^{4}}$$

by Lemma 4, we have

$$  \sum_{n=3}^\infty \, \frac{\log n }{n^{4}} \, = \, \int_3^\infty \, \frac{\log x}{x^4} \, d x \, + \, \alpha  $$

where $ \,  0 \, \leq \,  \alpha  \, \leq \, \frac{\log3}{3^4} $

 $$   \int_3^\infty \, \frac{\log x}{x^4} \, d\, x \, = \, - \, \frac{1}{3} \, \int_3^\infty \, \log x \,\, d \, x^{-3} \,
 = \, \frac{\log3}{3^4} \, + \, \frac{1}{3}\, \int_3^\infty \, x^{-4} \, d \, x \, $$

 $$ = \, \frac{\log3}{3^4} \, - \, \frac{1}{3^2} \, \int_3^\infty \, d  \, x^{-3} \, = \, \frac{\log3}{3^4} \, + \, \frac{1}{3^5} $$

therefore

 $$     \sum_{n=3}^\infty \, \frac{\log n }{n^{4}} \, \leq \, \frac{\log3}{3^4} \, + \, \frac{1}{3^5} \, + \, \frac{\log3}{3^4}  $$

 $$       \left| \, \, \zeta '( 4+ it ) \,\, \right| \, \geq \, \frac{\log2}{2^4} \, - \, \frac{2\log3}{3^4} \, -
 \, \frac{1}{3^5} \, \geq \, 0.012    $$

This completes the proof of Lemma 5.  \\

   Let $  \delta  =  \, \frac{1}{100} \, $,  $ \,\,  c_1 ,\, c_2 , \, ... \, , \,  $ is the  positive  constant. \\

{\bfseries\footnotesize LEMMA 6.}  When  $ \,\,  \sigma  \, \geq  \,  \frac{1}{2} , \,\,  |t| \geq  2  \,\, $, we have

$$  \left| \, \zeta ( \sigma  + it  )   \right|  \, \leq  \, c_1  \, |t|^\frac{1}{2}  $$

The lemma 6 follows from the  References [4], the theorem 2 of the page 140.      \\

{\bfseries\footnotesize LEMMA 7.} If $ \, f(z) \, $ is  the  analytic function  in the circle $ \,  | z - z_0  |\, \leq \, R  \,\,$,   $ \, 0 < r < R \,   $ , in the circle  $ \,\,   | z - z_0  |\, \leq \, r  \, $,  we have

$$  \left| \,  f(z) -  f(z_0)  \,  \right|  \, \leq  \,  \frac{2r}{R-r} \, \left ( \, A(R) - Re f(z_0) \,   \right)  $$

where  $ \,\,  A(R)\, = \, \max _{ | z - z_0  | \leq  R }   \,\,  Re f(z)  $.

 The lemma 7 follows from the  References [4], the theorem 2 of the page 61.     \\

   Now, we assume that Riemann hypothesis is correct, and abbreviation as RH. In other words, when $ \sigma >  \frac{1}{2}  $, the  function  $ \zeta ( \sigma  +  i t )  $  has no zeros. The function  $  \log  \zeta  ( \sigma + it )  $ is a multi-valued analytic function in  the  region $ \sigma >  \frac{1}{2} ,t \geq 1 $.  we choose  the principal branch of the function  $  \log  \zeta  ( \sigma + it )$,  therefore, if  $  \zeta  ( \sigma + it ) =1 ,  $  then  $  \log  \zeta  ( \sigma + it ) = 0 . $ \\

{\bfseries\footnotesize LEMMA 8.} If RH is  correct,   when $  \delta = \frac{1}{100} $,  $ \, \sigma\, \geq \, \frac{1}{2} \, + 2 \delta   \, ,\, \,  t  \geq  16 $, we have

   $$ \left | \,\, \log  \zeta ( \sigma  + it  ) \,\,  \right | \,\, \leq   \,\,  c_2    \log   t   +  c_3   $$

{\bfseries\footnotesize proof.} In Lemma 7, we choose  $  \,\, f(z) = \log \zeta ( z + 4 + it  ), \, z_0  = 0, \,   \,\, R = \frac{7}{2}
 - \delta  ,\,\, r = \frac{7}{2}  - 2 \delta , \,t \geq  16 \, $.  Because  $ \log \zeta ( z + 4 + it )  $  is the  analytic function  in the circle  $ | z | \leq  R  $,  by Lemma 7,   in the circle   $ | z | \leq  r   $, we have

$$  \left| \,\,  \log \zeta ( z + 4 + it  ) \,\, - \log \zeta (  4 + it  ) \,\,   \right| \,\, \leq \,\, \frac{7}{ \delta }\,\,
\left (  \,\, A(R)  - Re \log \zeta ( 4 + it ) \,\, \right)  $$

therefore

$$   \left| \,\,  \log \zeta ( z + 4 + it  ) \,\, \right|  \leq  \,\, \frac{7}{ \delta }  \,\, \left( \,\,  A(R)  +
 \left| \,\,  \log \zeta ( 4 + it )  \,\,  \right| \,\, \right)  +  \left| \,\, \log \zeta ( 4 + it ) \,\,  \right|  $$

by Lemma 6, we have

$$ A(R) = \max_{  | z - z_0 | \leq  R }   \log \left| \,\, \zeta ( z + 4 + it  ) \,\,  \right|  \leq \,\, \frac{1}{2} \, \log t
 \, + \, \log  c_1 $$

by  Lemma 5, we have

$$  \left| \,\,  \log \zeta ( z + 4 + it  ) \,\, \right| \,\, \leq  \,\, c_2  \, \log t \, +  \,  c_3  $$

therefore, when  $ \sigma  \geq  \frac{1}{2} + 2 \delta   $, we have

$$   \left| \,\,  \log \zeta ( \sigma + it  ) \,\, \right| \,\, \leq  \,\, c_2  \, \log t \, +  \,  c_3   $$

This completes the proof of Lemma 8. \\

{\bfseries\footnotesize LEMMA 9.} If  RH  is  correct,  when $ \delta = \frac{1}{100}  $,  $ \, t \geq  16 \, ,\, \rho = \frac{7}{2}  -  2 \delta $,  in the circle  $ |z| \leq  \rho   $, we have

$$  N \left( \,\,\rho \, ,\,\, \frac{1}{ \zeta ( z+ 4 + it  )\,\, - 1 } \,\, \right ) \,\, \leq \,\, \log\log t \, + \,  c_4  $$ \\

{\bfseries\footnotesize proof.} In Lemma 2,  we choose $ f(z) =  \,  \log \zeta ( z+ 4 + it  ) , \,\, R = \frac{7}{2} - \delta ,\,\, \rho = \frac{7}{2} - 2 \delta , \,\,a_\lambda \,\,\, (  \lambda = 1,2, ..., h ) $  are the zeros  of the function  $ \log \zeta ( z+ 4 + it  )   $  in the circle $  | z | < \rho $, each zero repeated according to their multiplicity.  Because the function $ \log \zeta ( z+ 4 + it  ) $   has no poles in the the circle  $  | z | < \rho $ ,  and  $ \log \zeta (  4 + it  )   $ is not equal to zero,  we have

$$ \log  \left| \,\log  \,  \zeta (  4 + it  ) \, \right|  = \, \frac{1}{ 2 \pi} \, \int_0^{ 2 \pi} \,  \log  \,
\left| \, \log\, \zeta (  4 + it + \rho e^{i \varphi} ) \, \right| \, d\varphi \, - \,  \sum_{\lambda = 1}^{h} \, \log
\, \frac{ \rho }{ |a_\lambda|  } $$

by Lemma 5 and Lemma 8, we have

$$  \sum_{\lambda = 1}^{h} \, \log  \, \frac{ \rho }{ |a_\lambda|  } \,\, \leq \,\,  \log\log t   \, + \, c_4  $$

because  $ z = 0  $   is neither  zero  nor  pole  of  the function $   \log \zeta ( z+ 4 + it  ) $,  if  $ r_0 $ is a sufficiently small   positive number,  we have

$$    \sum_{\lambda = 1}^{h} \, \log  \, \frac{ \rho }{ |a_\lambda| } \, = \, \int_{r_0}^\rho  \,
\left( \log \frac{\rho}{t} \right) \,\, d\, n(t, \frac{1}{f}  )  \,\, = \,\, \left[ \left( \log \frac{\rho}{t} \right) \,
   n(t, \frac{1}{f} ) \right] \bigg |_{r_0}^\rho $$

$$ + \,\, \int_{r_0}^\rho  \frac{n(t, \frac{1}{f}  )}{t } \,\, d \, t   \,\, = \,\,\int_{0}^\rho
 \frac{n(t, \frac{1}{f}  )}{t } \,\, d \, t  \,\, = \,\, N \left( \, \rho \, , \frac{1}{f} \,  \right )  $$

$$= \,\, N \left( \, \rho \, , \,\, \frac{1}{\log \zeta ( z+ 4 + it ) }  \, \right )  \,\, \geq \,\,  N \left( \, \rho \, , \,\,
\frac{1}{ \zeta ( z+ 4 + it ) - 1  }  \, \right )  $$

This completes the proof of Lemma 9. \\

{\bfseries\footnotesize THEOREM.} \,  If  RH  is  correct,   when $  \sigma  \, \geq \, \frac{1}{2}  \,  +  \, 4 \delta \, ,\,\,  \delta =
\frac{1}{100}\,, \,\, t \geq  16    \,\,  $, we have

$$  \left| \,  \zeta ( \sigma + it  ) \, \right| \,  \leq    c_{8}  \,\left( \, \log  t  \, \right)^{c_{6}}  $$ \\

{\bfseries\footnotesize proof.}  In Lemma 3,  we choose $ f(z) \,= \, \zeta ( z + 4 + it )  , \, t \geq 16  $,  $ R \, = \,  \frac{7}{2} \, - \, 2 \delta , \,\, r \, = \, \frac{7}{2} \,  - \, 3 \delta  $.  by  Lemma 5,  we have  $\, f(0) \, = \, \zeta ( 4 + it ) \, \neq \,  0, \, \infty , \,  1, \,\,$ and  $  | f'(0)| \, = \, | \zeta ' (  4 + it) |\, \geq  0.012  \, , \,\,\,\,  | \, f(0) \, |  \, = \, \left| \, \zeta ( 4 + it ) \, \right|  \, \leq  \, 1.0824 $. because $ \zeta ( z + 4 +it )  $  is  the analytic function,  and it have neither  zeros nor poles in the circle $ |z| \leq  R $,   we have

$$   N \left( \, R \, , \, \frac{1}{f}  \right)  \,  =\, 0 \,\, , \,\,\,\,\,\,\,  \,\,\,\,  N  \left( \, R \, , \, f  \right )  =  0    $$

therefore, by Lemma 9, we have

$$  T \, \left( \, r \, , \zeta ( z + 4 + it  )  \, \right ) \,\,  \leq  \,\,  2 \log\log t \, + \, c_5    $$

In Lemma 1,  we choose $ \, R = \frac{7}{2} - 2\delta \,, \,  \rho \, = \, \frac{7}{2}  \, - \,  3 \delta  , \,\,  r \, = \,  \frac{7}{2} \, -
\, 4 \delta   $.  by the maximal principle,  in the circle  $ |z | \leq  r   $,  we have

$$  \log^+  \left| \,\zeta ( z + 4 + it  ) \,  \right|   \, \leq  \, c_6 \, \log\log t \, + \, c_7   $$

therefore,  when  $ \sigma  \geq \frac{1}{2} \, + \, 4 \delta   $, we have

$$   \log^+  \left| \,\zeta ( \sigma + it  ) \,  \right|   \, \leq  \, c_6 \, \log\log t \, + \, c_7   $$

$$  \log \,  \left| \,\zeta ( \sigma + it  ) \,  \right|   \, \leq  \, c_6 \, \log\log t \, + \, c_7     $$

$$  \left| \,  \zeta ( \sigma + it  ) \, \right| \,  \leq    c_{8}  \,\left( \, \log  t \, \right)^{c_{6}}  $$

This completes the proof of  Theorem. \\

\vspace{20mm} \centerline{ REFERENCES } \vspace{5mm}

[1]  Zhuang Q.T, { \itshape Singular   direction   of   meromorphic   function,}  BeiJing, Science Press, 1982. ( in Chinese )\\

[2]  Yang Lo, { \itshape Value distribution theory, } Springer, 1993. \\

[3]  Hua  L.G, { \itshape  Introduction of number theory,}   BeiJing,  Science Press, 1979. ( in Chinese ) \\

[4]  Pan C.D, Pan C.B,{ \itshape Fundamentals of analytic number theory,  }  BeiJing,  Science Press,  1999. ( in Chinese )  \\

[5]  Zhuang  Q.T and   Zhang  N.Y, { \itshape  Complex variables functions,}   BeiJing,  Peking University press ,   1984. ( in Chinese )  \\

[6]  Hua  L.G, { \itshape  Introduction of advanced mathematics ( first book  of the second volume   ), }   BeiJing,  Science Press, 1981. ( in Chinese )\\

[7] William Cherry and Zhuan Ye,  { \itshape  Nevanlinna's Theory of Value Distribution, }  Springer,  2001.\\

[8] Serge Lang and William Cherry, { \itshape   Topics in Nevanlinna Theory, }  Springer,  1990.\\

[9] Matthew M.Buck,  { \itshape  Nevanlinna Theory, } Final Report,   2007.

\end{document}